\documentclass[10pt]{article}
\usepackage{fullpage,amsfonts,amsmath,amssymb,amsopn,amsthm,mdwlist,wasysym,url,graphicx,verbatim}
\usepackage[cns,uns]{llconnex}

\input diagxy

\DeclareMathOperator{\rep}{rep}

\newcommand{\op}{{\rm op}}

\newcommand{\Aut}{{\bf Aut}}
\newcommand{\Cat}{{\bf Cat}}
\newcommand{\Vect}{{\bf Vec}}

\newcommand{\Mon}{{\bf Mon}}

\newcommand{\ladj}{\dashv}

\newcommand{\too}{\longrightarrow}

\newcommand{\oto}[1]{\stackrel{\underrightarrow{\stackrel{#1}{\,\,\,\,\hphantom{\too}}}}{}}  
 
\author{Micah Blake McCurdy}
\date{March 2, 2009}
\title{Tannaka Reconstruction of Weak Hopf Algebras in Arbitrary Monoidal Categories}

\begin{document}

\maketitle

\begin{abstract}
\noindent We introduce a variant on the graphical calculus of Cockett and Seely\cite{CS} for monoidal functors and illustrate it with a 
discussion of Tannaka 
reconstruction, some of which is known and some of which is new. The new portion is: given a separable Frobenius functor $F: A \too B$ from a
monoidal category $A$ to a suitably complete or cocomplete braided autonomous category $B$, the usual formula for Tannaka reconstruction gives a 
weak bialgebra in $B$; if, moreover, $A$ is autonomous, this weak bialgebra is in fact a weak Hopf algebra.
\end{abstract}

\section{Introduction}
Broadly speaking, Tannaka duality describes the relationship between algebraic objects and their representations, for an excellent introduction, 
see \cite{JoyalStreet}. On the one hand, given an algebraic object $H$ in a monoidal category $B$ (for instance, a hopf object in the category
$\Vect_k$ of vector spaces over a field $k$), one can consider the functor which takes algebraic objects of the given type to their category of
representations, $\rep_B H$, for which there is a canonical forgetful functor back to $B$. This process is \emph{representation} and it can be defined in a
great variety of situations, with very mild assumptions on $B$. For instance:
\begin{itemize}
	\item If $J$ is a bialgebra in a braided linearly distributive category $B$, then $\rep_B J$ is linearly distributive \cite{CS}.
	\item If $H$ is a Hopf algebra in a braided star-autonomous category $B$, then $\rep_B H$ is star-autonomous \cite{CS}.
\end{itemize} Note that the first example includes the familiar notion of a bialgebra in a monoidal category giving rise to a monoidal category
of representations, since every monoidal category is degenerately linearly distributive by taking both of the monoidal products to be the same.

On the other hand, given a suitable functor $F: A \too B$, we can try to use the properties of $F$ (which of course include those of $A$ and $B$)
to build an algebraic object in $B$; this is called (Tannaka) \emph{reconstruction}, since historically the algebraic objects have been considered 
primitive. We denote the reconstructed object as $E_F$, following \cite{Street}.
This requires making more stringent assumption on $B$; certainly it must be braided; it must be autonomous, and it is assumed that $B$ admits 
certain ends or coends which cohere with the monoidal product. For instance, under these assumptions: \begin{itemize}
	\item If $F$ is any functor, then $E_F$ is a monoid in $B$.
	\item If $A$ is monoidal and $F$ is a monoidal functor, then $E_F$ is a bialgebra in $B$.
	\item If $A$ is autonomous and $F$ is a strong monoidal functor, then $E_F$ is a Hopf algebra in $B$ (\cite{Schauenberg}, \cite{Majid}).
\end{itemize}
To this we add the following: \begin{itemize}
	\item If $A$ is monoidal and $F$ is a separable Frobenius functor, then $E_F$ is a weak bialgebra in $B$.
	\item If $A$ is autonomous and $F$ is a separable Frobenius functor, then $E_F$ is a weak Hopf algebra in $B$.
\end{itemize}

Note that the notion of weak Hopf algebra considered by Haring-Oldberg \cite{HO} is a different notion, the sense of weak we use here is that
of \cite{BNS}. It should be noted that a special case of the weak Hopf algebra result has been obtained by
Pfeiffer\cite{P}, where $A$ is taken to be a modular category and $B$ is taken to be $\Vect_k$; see also \cite{Sz}.

In favourable circumstances, reconstruction is left adjoint to representation, for instance:

\[ \bfig
\node a(0,0)[{\rm Monoids}(B)]
\node b(0,-1000)[\Cat/B]
\arrow|r|/{@{>}@/^1em/}/[a`b;\rep_B-]
\arrow|l|/{@{>}@/^1em/}/[b`a;E_-]
\place(0,-500)[\ladj]

\node a(+1500,0)[{\rm Bialgebras}(B)]
\node b(+1500,-1000)[\Mon\Cat_{\rm Strong}/B]
\arrow|r|/{@{>}@/^1em/}/[a`b;\rep_B-]
\arrow|l|/{@{>}@/^1em/}/[b`a;E_-]
\place(+1500,-500)[\ladj]

\node a(+3000,0)[{\rm Hopf}(B)]
\node b(+3000,-1000)[\Aut\Cat_{\rm Strong}/B]
\arrow|r|/{@{>}@/^1em/}/[a`b;\rep_B-]
\arrow|l|/{@{>}@/^1em/}/[b`a;E_-]
\place(+3000,-500)[\ladj]

\efig \]

In the weak cases which I discuss in the sequel, however, we do not have these adjunctions.

\section{Graphical Notation for Monoidal and Comonoidal Functors}
Before we discuss the reconstruction itself, we discuss notations for monoidal and comonoidal functors. The original notion for graphically 
depicting monoidal functors
as transparent boxes in string diagrams is due to Cockett and Seely\cite{CS}, and has recently been revived and popularized by 
Mellies\cite{Mellies} with prettier graphics and an
excellent pair of example calculations which nicely show the worth of the notation. However, a small modification improves the notation considerably.
For a monoidal functor $F : A \too B$, we have a pair of maps, $Fx \tens Fy \too F(x \tens y)$ and $\e \too F\e$, which we notate as follows:
\begin{center}\resizebox{!}{50pt}{\includegraphics{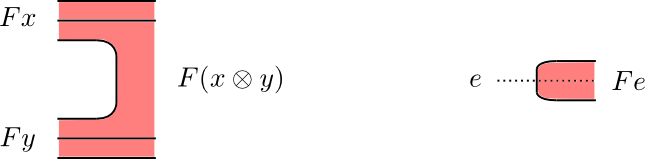}}\end{center} 
Similarly, for a comonoidal $F$, we have maps $F(x \tens y) \too Fx \tens Fy$ and $F\e \too e$ which we notate in the obvious dual way, as follows:
\begin{center}\resizebox{!}{50pt}{\includegraphics{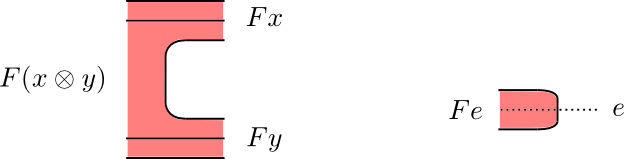}}\end{center} 

Graphically, the axioms for a monoidal functor are depicted as follows:
\begin{center}\resizebox{450pt}{!}{\includegraphics{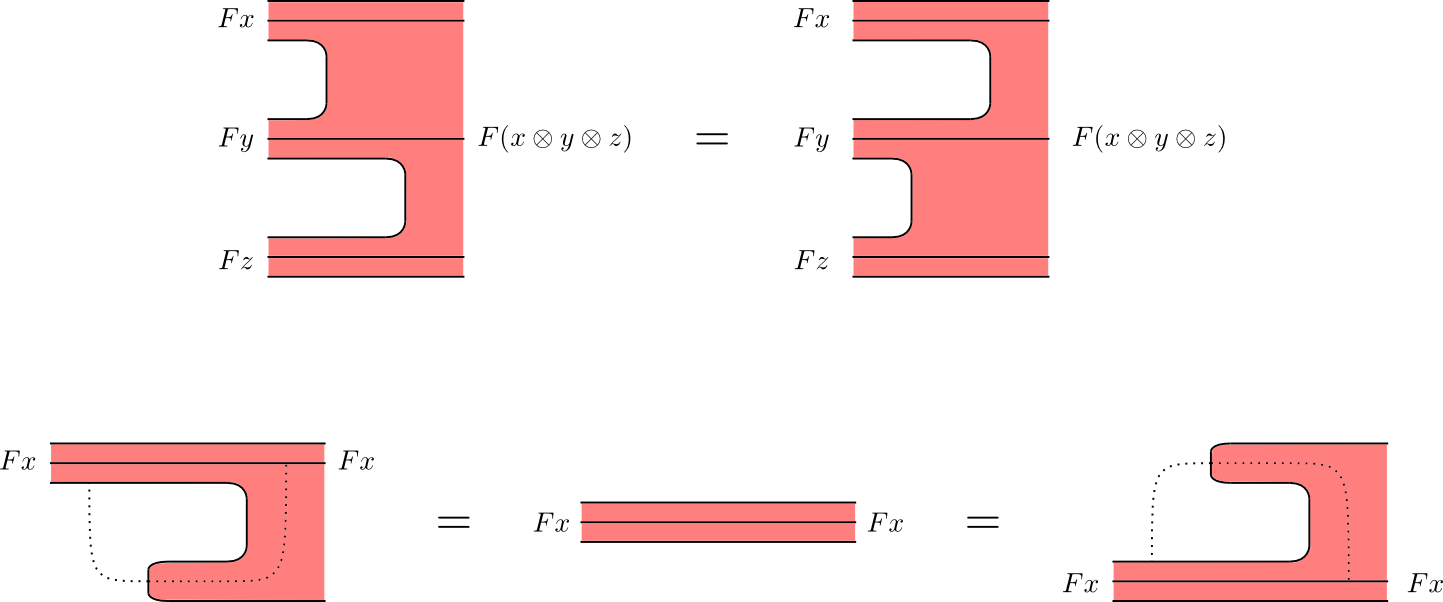}}\end{center}
where, once again, the similar constraints for a comonoidal functor are exactly the above with composition read right-to-left instead of left-to-right.

It is curious and pleasing that the two unit axioms bear a superficial resemblence to triangle-identies being applied along the 
boundary of the ``$F$-region''.

The above axioms seem to indicate some sort of ``invariance under continuous deformation of $F$-regions''. For a functor which is both 
monoidal and comonoidal, pursuing this line of thinking leads one to consider the following pair of axioms:
\begin{center}\resizebox{450pt}{!}{\includegraphics{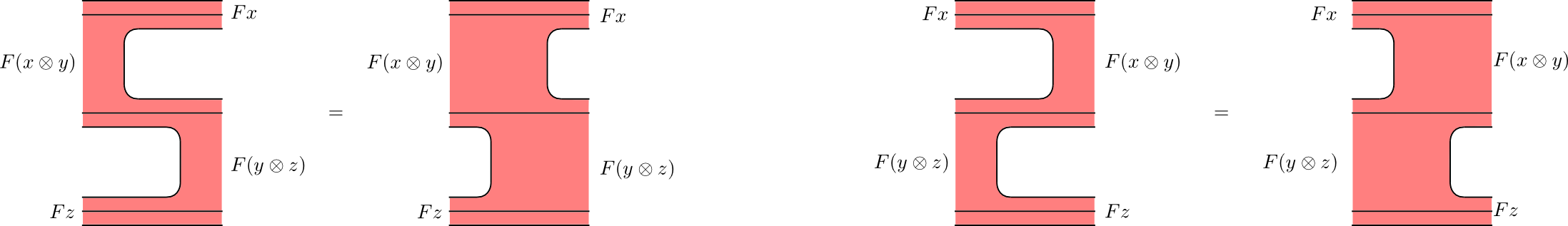}}\end{center}
A functor with these properties has been called a Frobenius monoidal functor by Day and Pastro \cite{DP}; moreover, it is the notion which
results from considering linear functors in the sense of \cite{CS} between degenerate linearly distributive categories. They are more common than strong
monoidal functors yet still share the key property of preserving duals. 

Let $F: A \too B$ be a (mere) functor between monoidal categories, where $B$ is assumed to be left closed. Then define
\[ E_F = \int_{a \in A} \left[ Fa, Fa \right] \] where I assume that $A$ and $B$ are such that the indicated end exists. As Richard Garner imimitably asked at the 2006 PSSL
in Nice, ``Have you considered enriching everything?''. I do not discuss the matter here, but it has been considered by Brian Day\cite{Day}.

There is a canonical action of $E_F$ on $Fx$ for each object $x$ in $A$, which we denote as $\alpha = \alpha_x : E_F \tens Fx \too Fx$.
This is defined as:
\[ E_F \tens Fx = \left( \int_{a \in A} [Fa,Fa] \right) \tens Fx \oto{\pi_x \tens Fx} [Fx,Fx] \tens Fx \oto{\rm ev} Fx \]
using the $x$'th projection from the end followed by the evaluation of the monoidal closed structure of $B$. 
The dinaturality of the end in $a$ gives rise to
the naturality of the action on $Fa$ in $a$, which we notate as:
\begin{center}\resizebox{300pt}{!}{\includegraphics{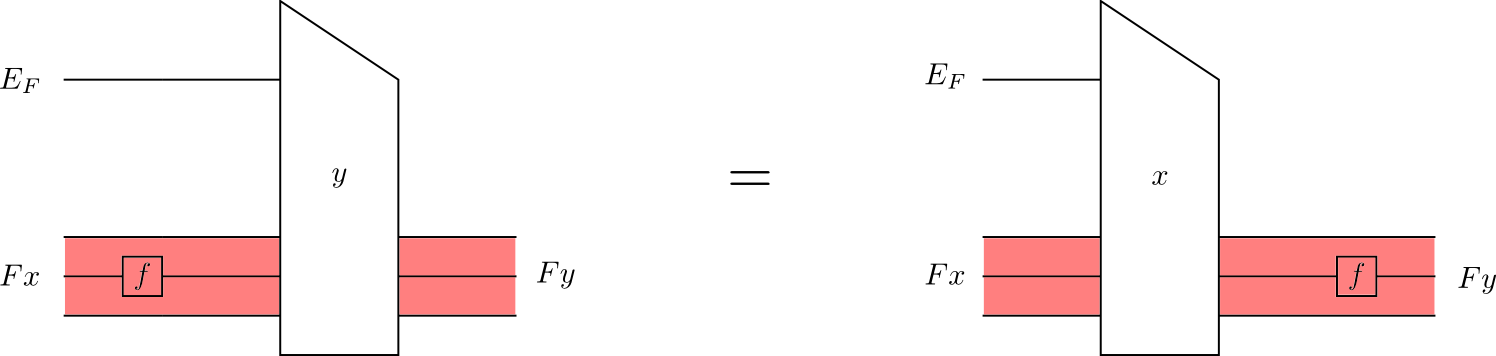}}\end{center}
Let us now assume that the closed structure of $B$ is given by left duals, that is, $[a,b] = b \tens La$.
If we also assume that $B$ is braided and that the tensor product coheres with the ends in $B$, then we obtain canonical actions of 
$E_F^n$ on $Fx_1 \otimes \ldots \otimes Fx_n$, written $\alpha^n_x$. Taking $\alpha^1 = \alpha$, we define $\alpha^n$ recursively 
as follows:
\[ E_F^n \tens \left( \bigotimes_{j=1}^n Fx_j \right) =
E_F^{n-1} \tens E_F \tens \left( \bigotimes_{j=1}^{n-1} Fx_j \right) \tens Fx_n \oto{\rm braid}
E_F^{n-1} \tens \left( \bigotimes_{j=1}^{n-1} Fx_j \right) \tens E_F \tens Fx_n 
\oto{\alpha^{n-1} \tens \alpha^1} \bigotimes_{i=1}^n Fx_i \]

For any map $f: X \too E_F^n$, we may paste together $f$ with $\alpha^n$ along $\bigotimes_{i=1}^n Fx_i$ to obtain a ``discharged form'' of 
$f$: \[ X \tens Fx_1 \tens \ldots \tens Fx_n \too Fx_1 \tens Fx_1 \tens \cdots \tens Fx_n \]
Two maps are equal if and only if they have the same discharged forms.

Many treatments instead consider \[ E^F = \int^{a \in A} [Fa,Fa] \] Under dual assumptions to the ones above -- namely, that the coend exists and
coheres with the tensor product in $B$ -- we obtain canonical coactions of $E^F$ on $Fx$ for each $x \in A$, and iterated coactions, \&c. This
approach has certain technical benefits; among others, that the tensor product in $B = \rm Vect$ coheres with coends but not with ends. However,
the notation we use covers \emph{both} cases, $E_F$ and $E^F$. For the former, one must read composition left-to-right, and for the latter, from
right-to-left. We write $E_F$ as a label for convenience, preferring for convenience to read in the conventional English way, but it is a crucial 
feature of the notation that in fact no choice is made.

Without assuming that $F$ bears a monoidal structure, one can define a monoid structure on $E_F$, as follows:
\begin{center}\resizebox{450pt}{!}{\includegraphics{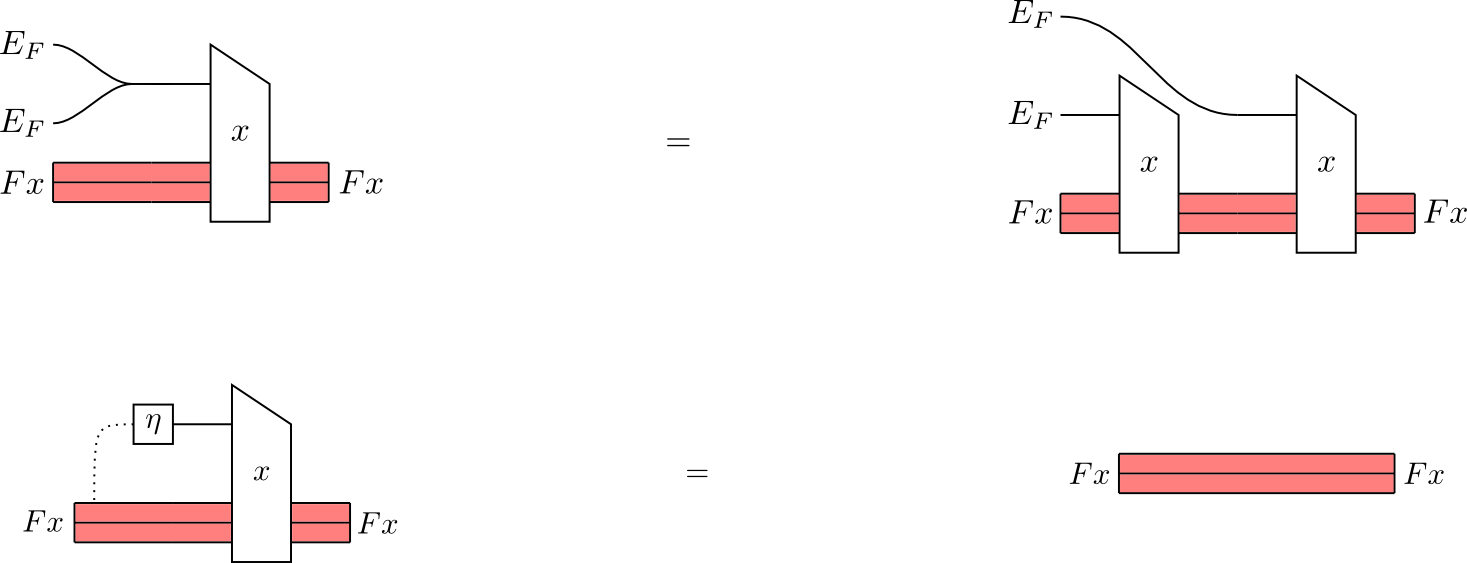}}\end{center} 
Note that this monoidal structure is associative and unital, without any assumption on $F$.

Furthermore, if $F$ is known to be (lax) monoidal and comonoidal
(without at the moment assuming any coherence between these structures) we can define a comonoid structure on $E_F$. 
\begin{center}\resizebox{450pt}{!}{\includegraphics{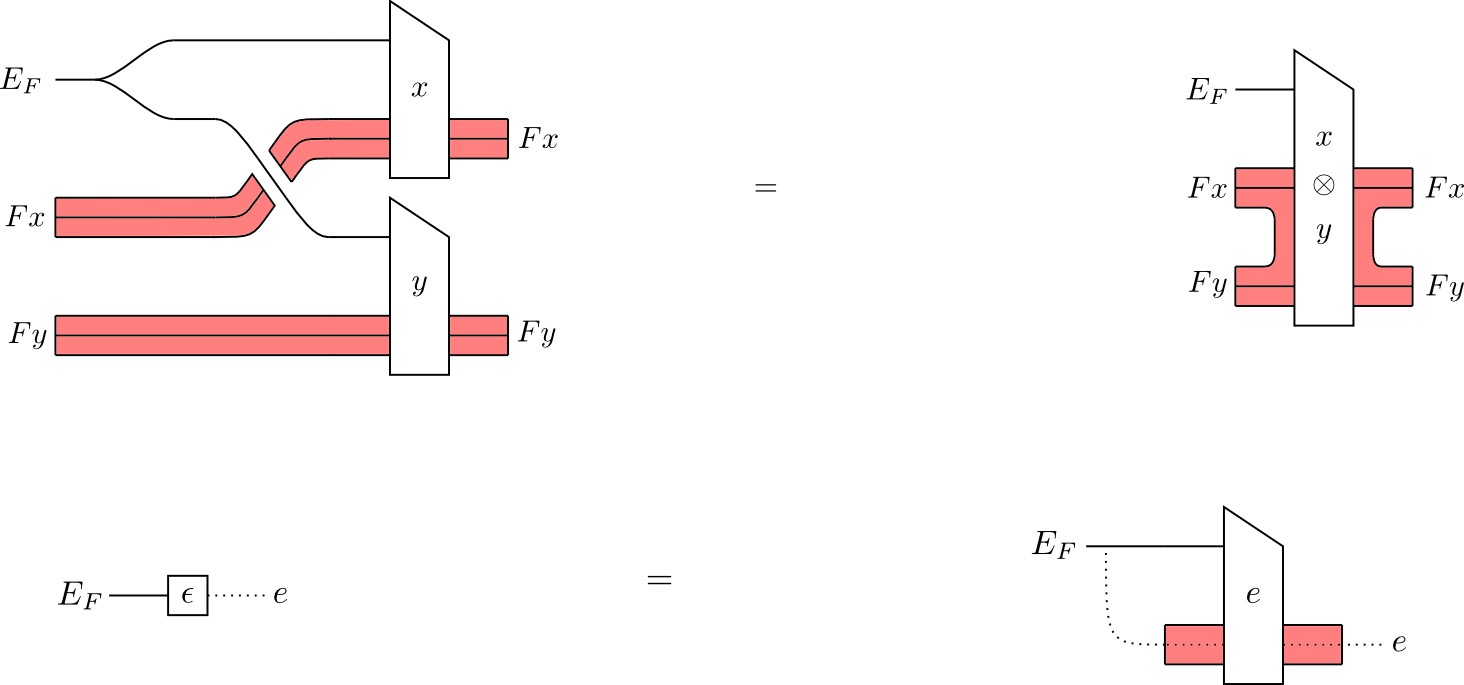}}\end{center} 
Finally, if $A$ is known to have (left, say)
duals, we can define a canonical map $S: E_F \too E_F$ which we think of as a candidate for an antipode.
\begin{center}\resizebox{450pt}{!}{\includegraphics{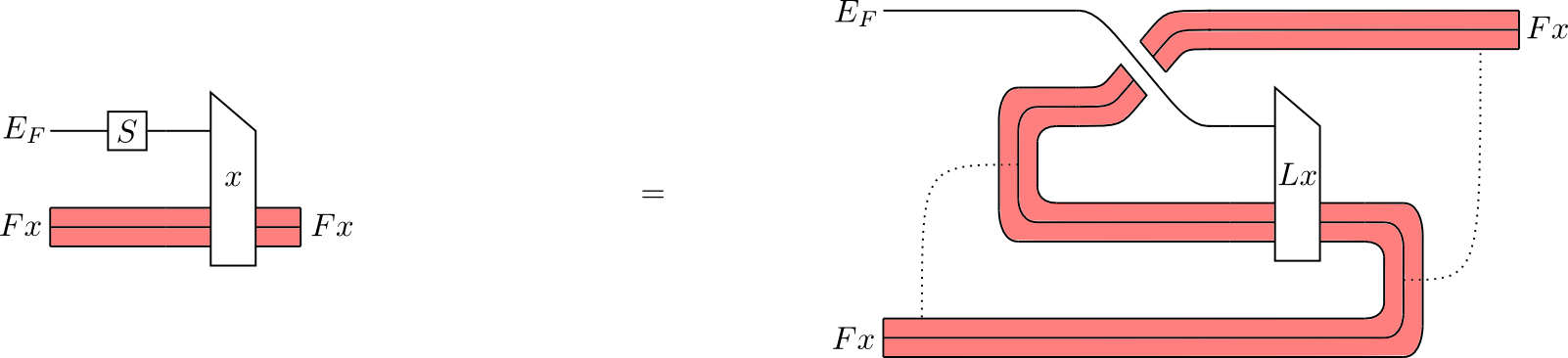}}\end{center}

Notice in particular how the monoidal and comonoidal structures on $F$ permit one to consider the application of $F$ as not merely ``boxes'' but more like a
flexible sheath.

Now, the above is the raw data for two different structures, namely, Hopf algebras and weak Hopf algebras, which differ only in axioms. 
It has been remarked
before that requiring $F$ to be strong (that is, in our treatment, demanding that the monoidal and comnoidal structures be mutually inverse) 
makes the above data into a Hopf algebra. Before we discuss the Hopf data (that is, the antipode), let us first consider the bialgebra data.
A bialgebra in a braided monoidal category satisfies the following four axioms. First, the unit followed by the counit must be the identity:
\begin{center}\resizebox{300pt}{!}{\includegraphics{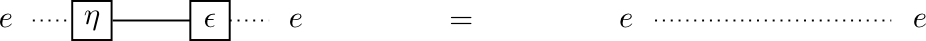}}\end{center}
Second and third, the unit and counit must respect the comultiplication and multiplication, respectively:
\begin{center}\resizebox{450pt}{!}{\includegraphics{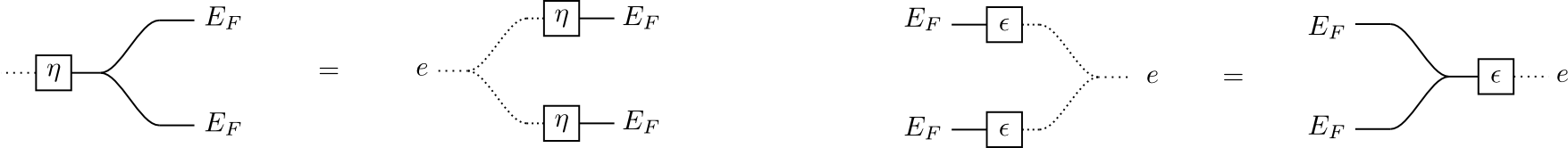}}\end{center}
Fourthly, the multiplication must cohere with the comultiplication, with the help of the braiding:
\begin{center}\resizebox{300pt}{!}{\includegraphics{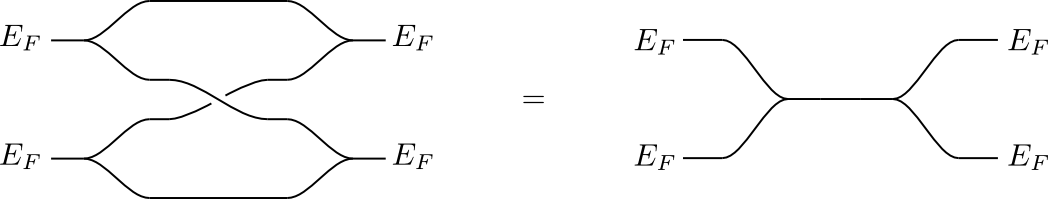}}\end{center}

Some easy calculations show how the strength of $F$ features crucially in showing all four of these axioms.

For the first of these, we calculate:
\begin{center}\resizebox{450pt}{!}{\includegraphics{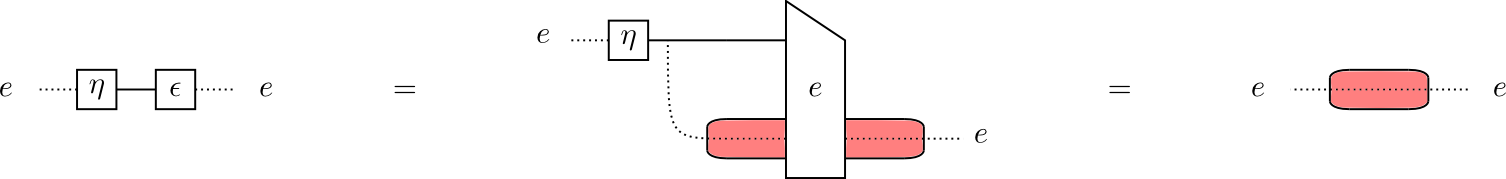}}\end{center}
and we see that this composite is the identity on $e$ precisely when $e \too Fe \too e$ is the identity.

For the second bialgebra axiom, we have the following two calculations:
\begin{center}\resizebox{450pt}{!}{\includegraphics{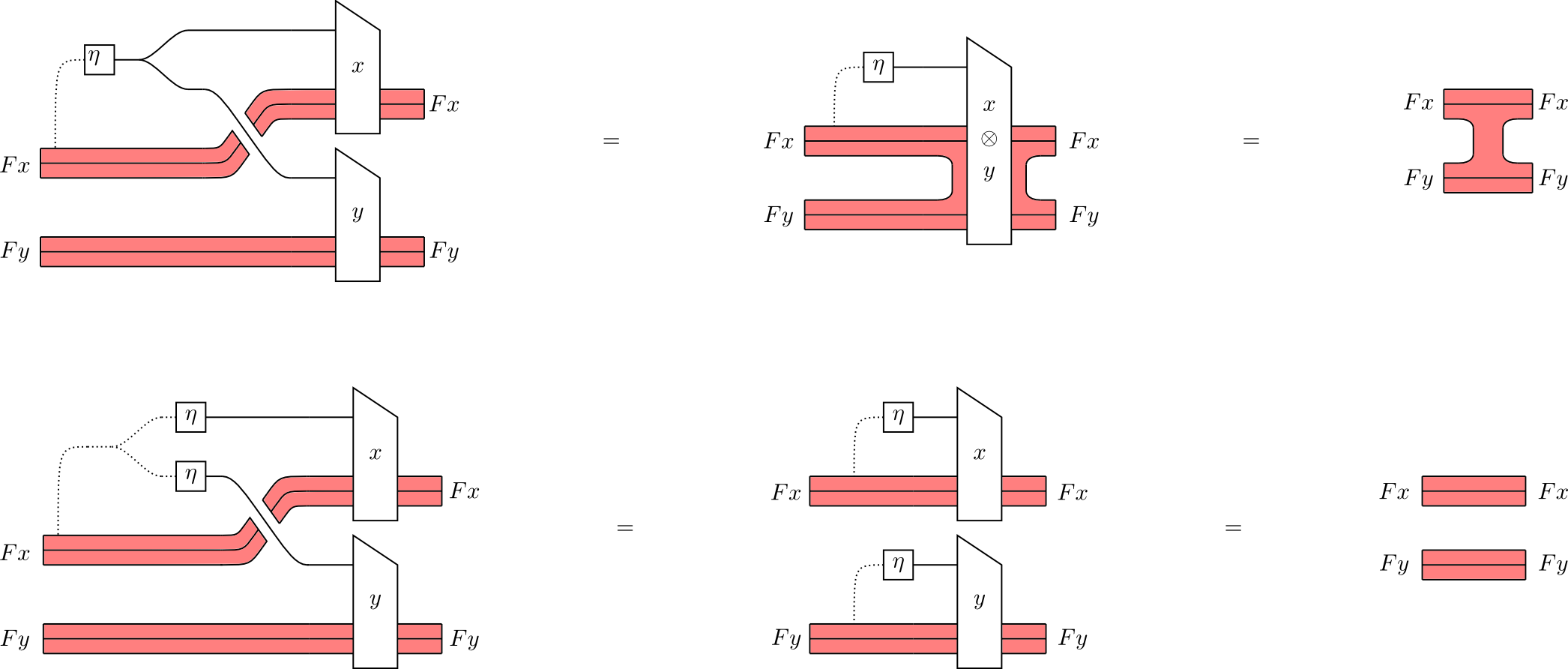}}\end{center}
and so we see that these two are equal precisely when $Fx \tens Fy \too F(x \tens y) \too Fx \tens Fy$ is the identity.

For the third bialgebra axiom, we have the following two calculations:
\begin{center}\resizebox{450pt}{!}{\includegraphics{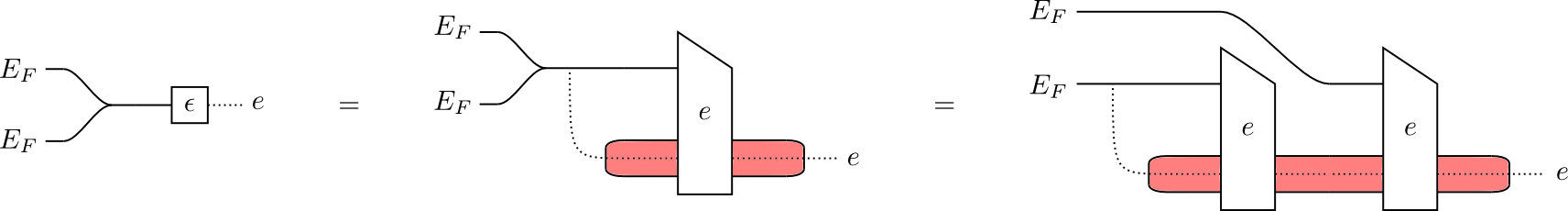}}\end{center}
\begin{center}\resizebox{450pt}{!}{\includegraphics{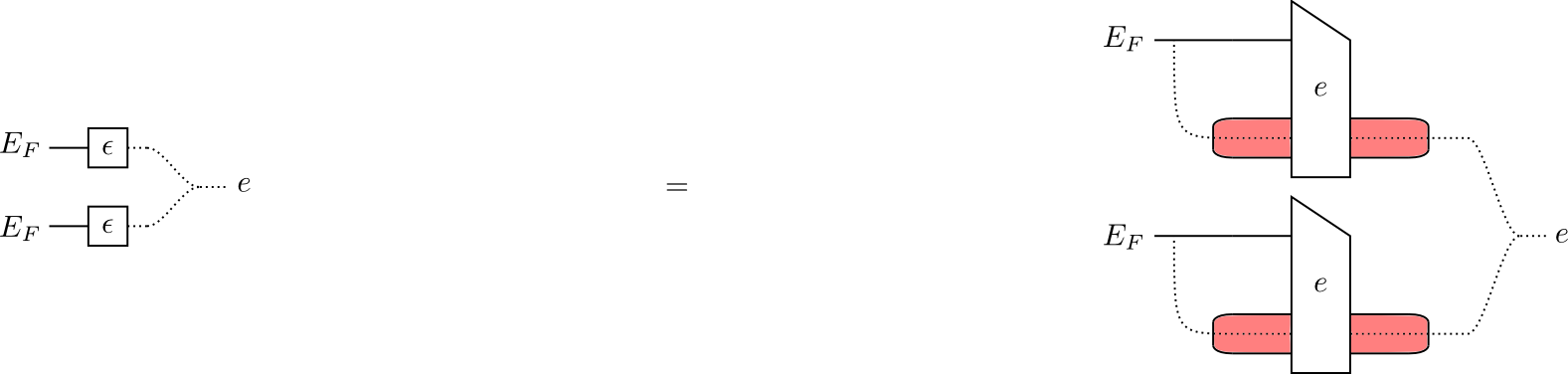}}\end{center}
and we see that for these two to be equal, it suffices to have $Fe \too e \too Fe$ be the identity, the use of which between the two actions in the first calculation
gives the result.

\begin{figure}[ht]
\begin{center}\resizebox{!}{550pt}{\includegraphics{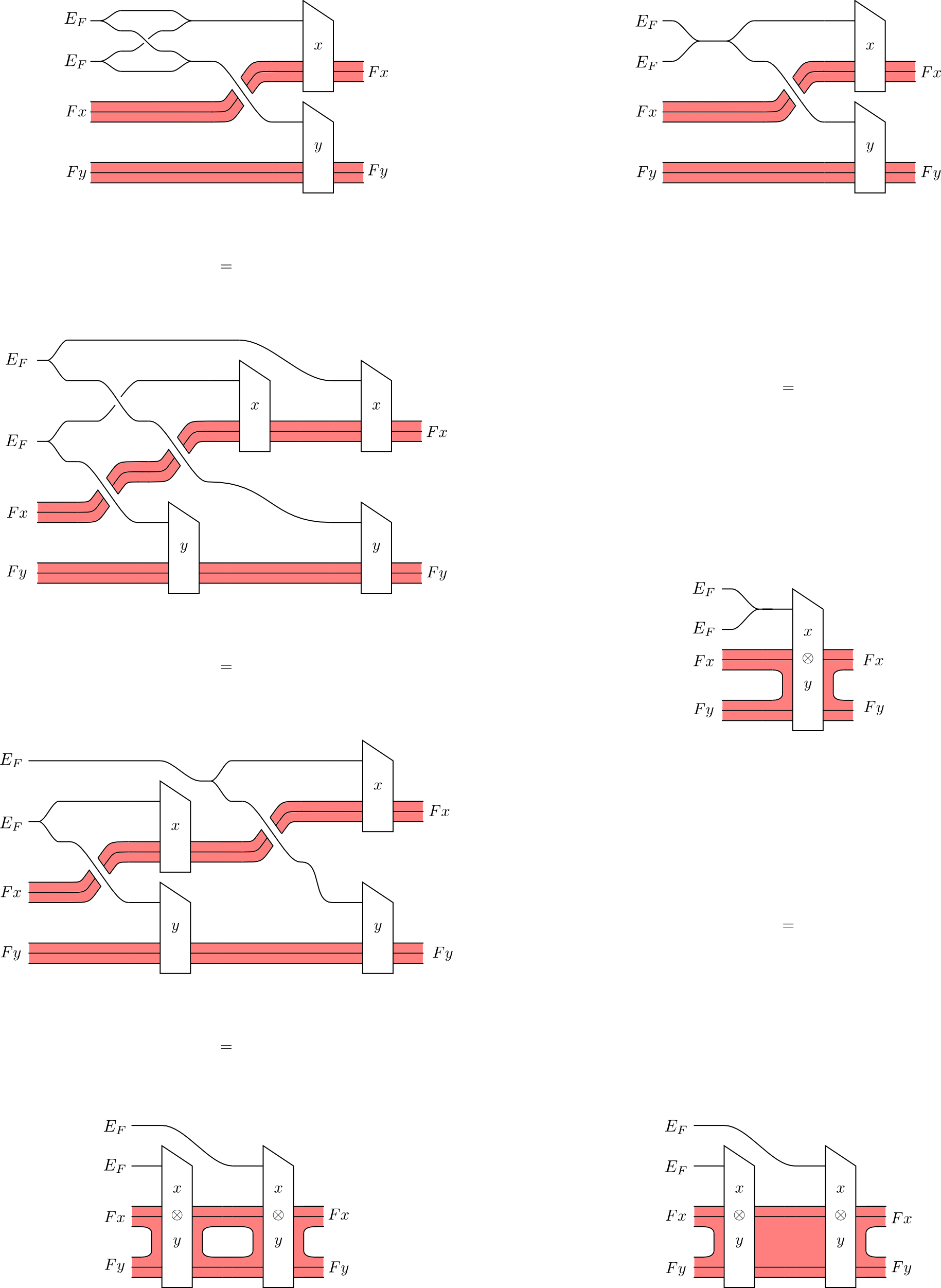}}\end{center}
\caption{Coherence of the multiplication with the comultiplication}
\label{bialgebra1}
\end{figure}

Finally, for the final bialgebra axiom, the calculations shown in figure (\ref{bialgebra1}) compute the discharged forms as
\begin{center}\resizebox{350pt}{!}{\includegraphics{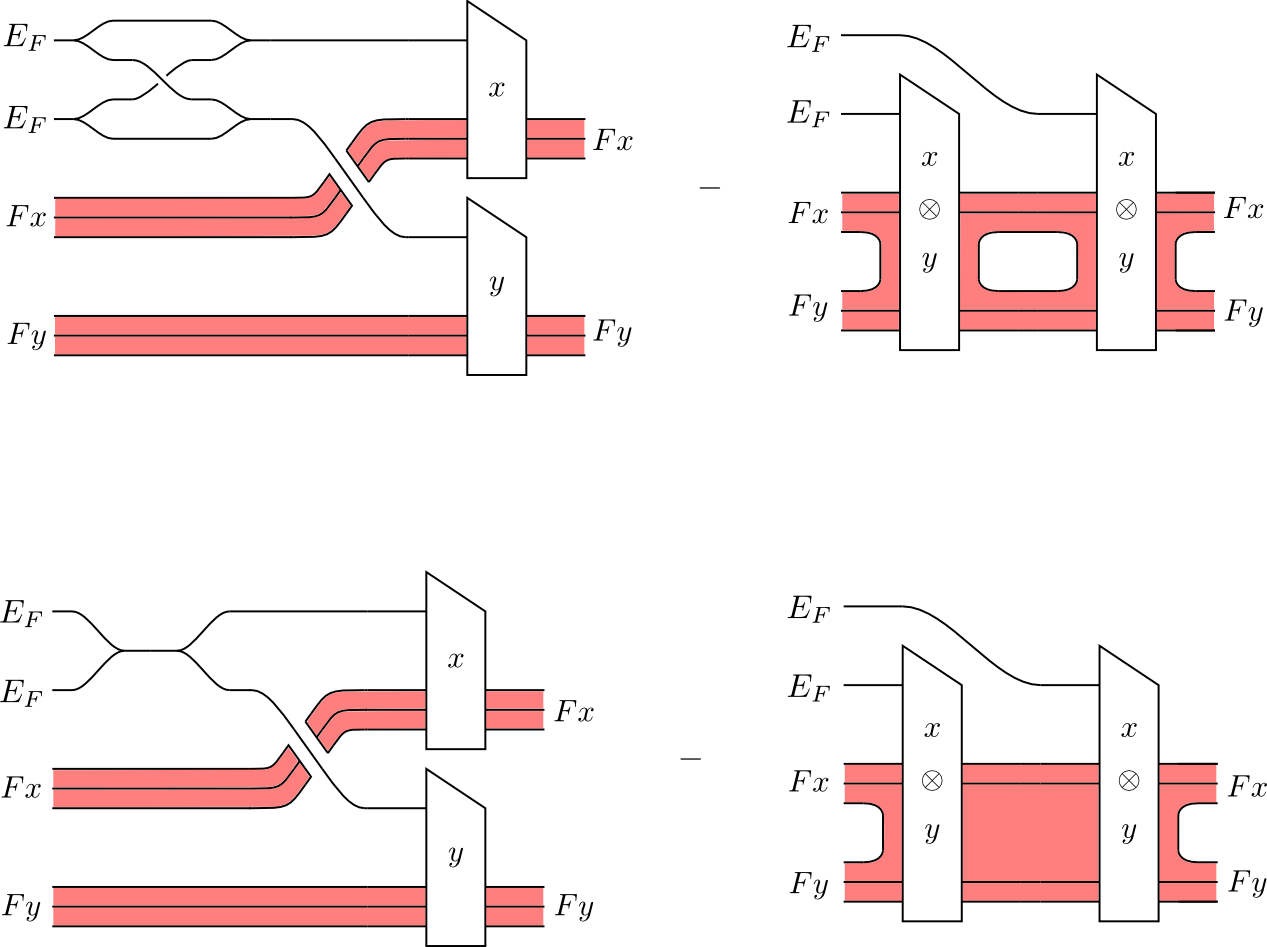}}\end{center}
which shows that it suffices to request that $F(x \tens y) \too Fx \tens Fy \too F(x \tens y)$ should be the identity.

Demanding that $F$ be strong imposes the following four conditions on the monoidal/comonoidal structure:
\begin{center}\resizebox{450pt}{!}{\includegraphics{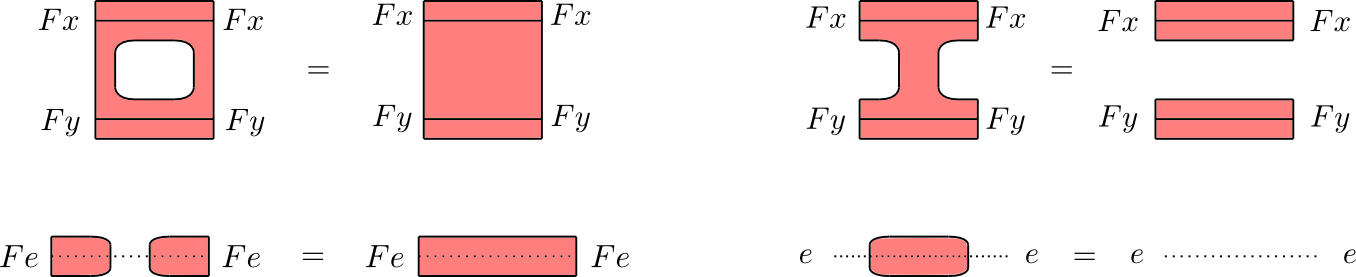}}\end{center}

Notice that precisely one of them preserves the number of connected components of $F$, namely, the one which is used in the proof of the coherence of the multiplication with
the comultiplication. A monoidal functor satisfying this axiom has been called ``separable'' by some.
To move from a (strong) bialgebra to a weak one, this axiom is the only one which is retained. The coherence of the unit with the counit is discarded entirely, and the second
and third axioms are replaced with the following four axioms:
\begin{center}\resizebox{450pt}{!}{\includegraphics{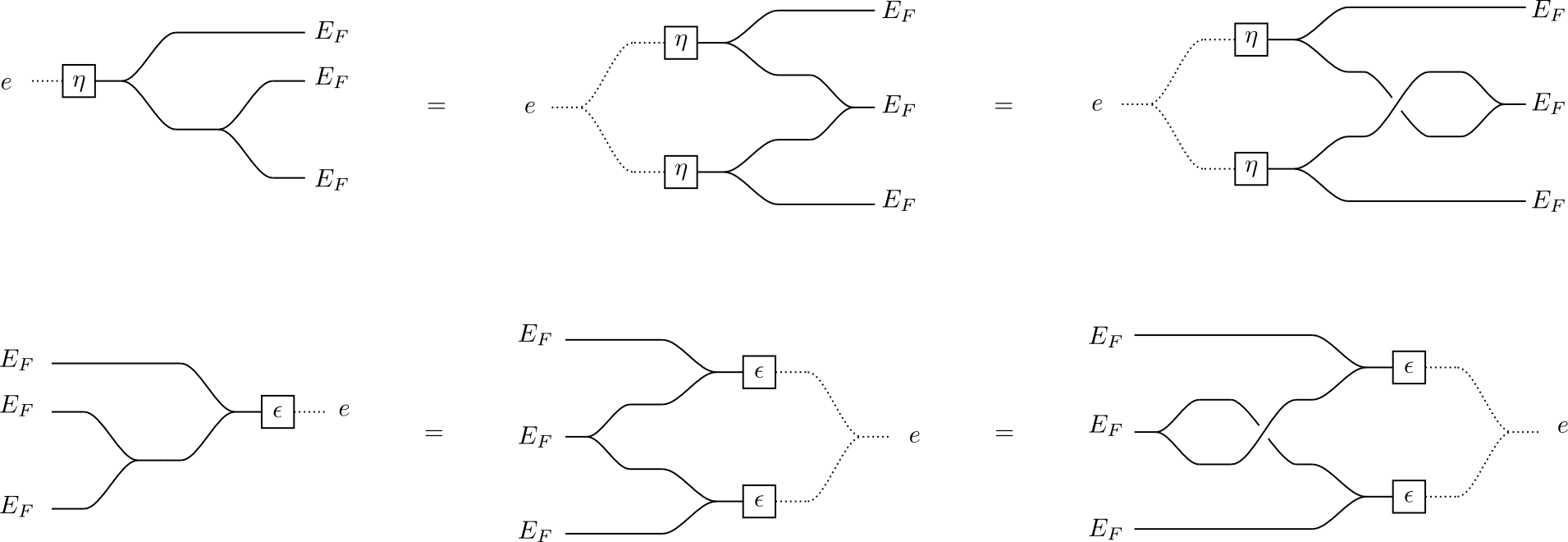}}\end{center}

We first examine the unit axioms. In discharged form, the first unit expression is calculated as:
\begin{center}\resizebox{450pt}{!}{\includegraphics{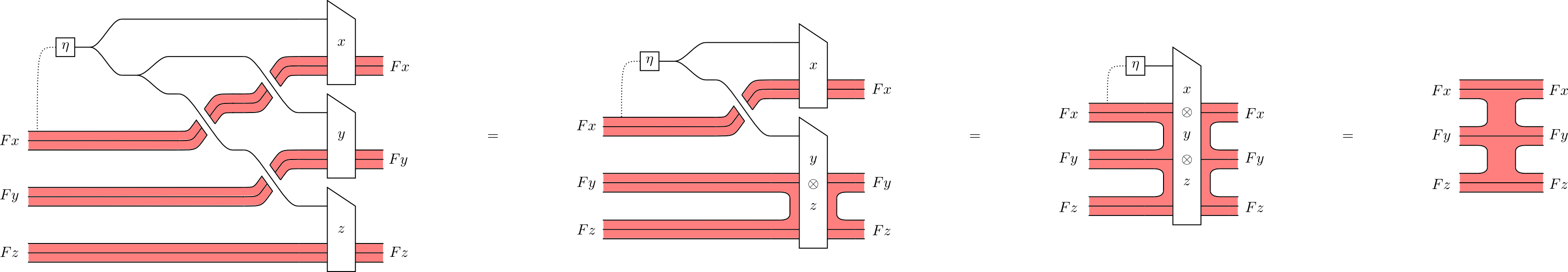}}\end{center}

\begin{figure}[ht]
\begin{center}\resizebox{!}{550pt}{\includegraphics{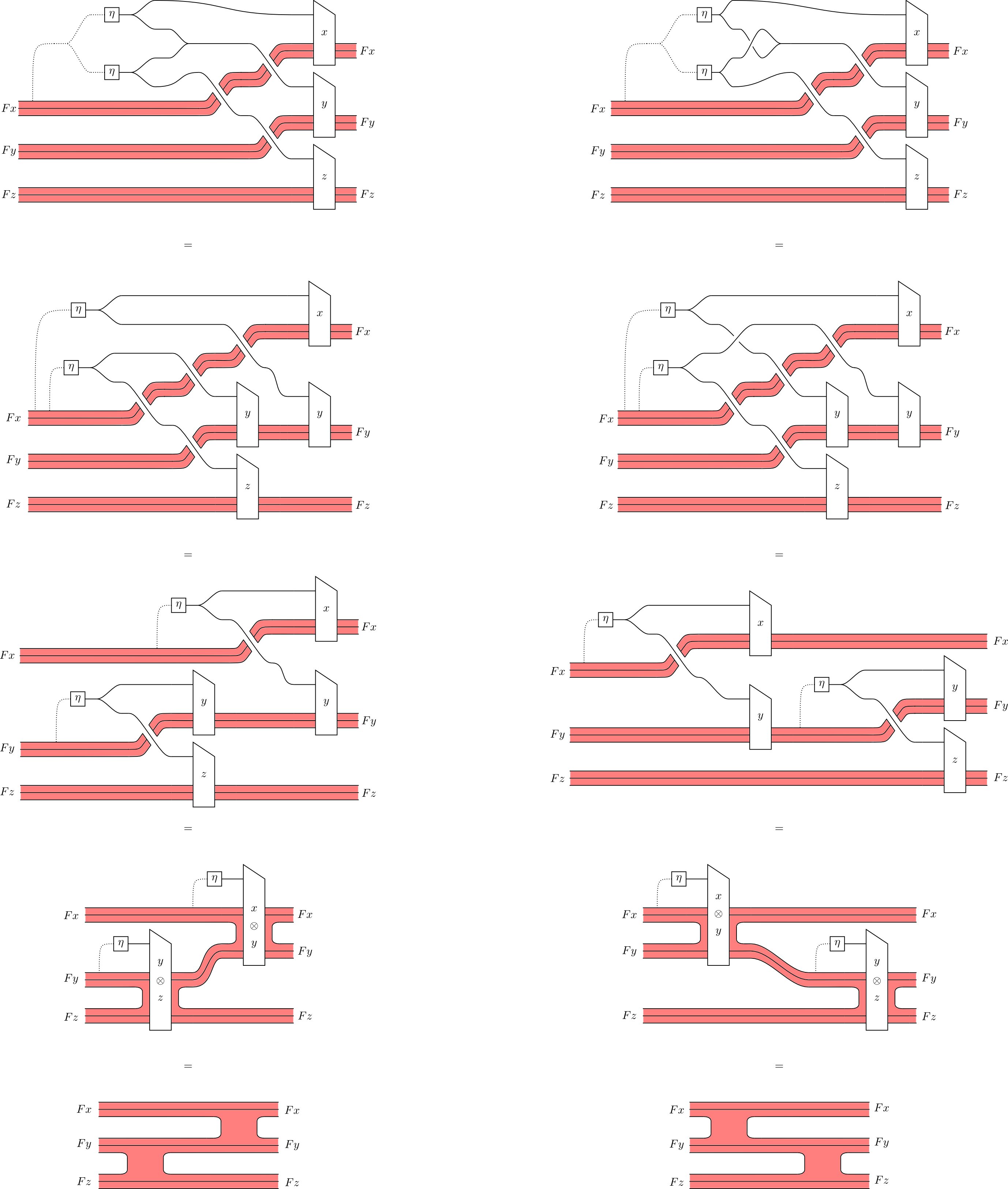}}\end{center}
\caption{Weak unit calculations}
\label{weakunit23}
\end{figure}

\begin{figure}[ht]
\begin{center}\resizebox{!}{550pt}{\includegraphics{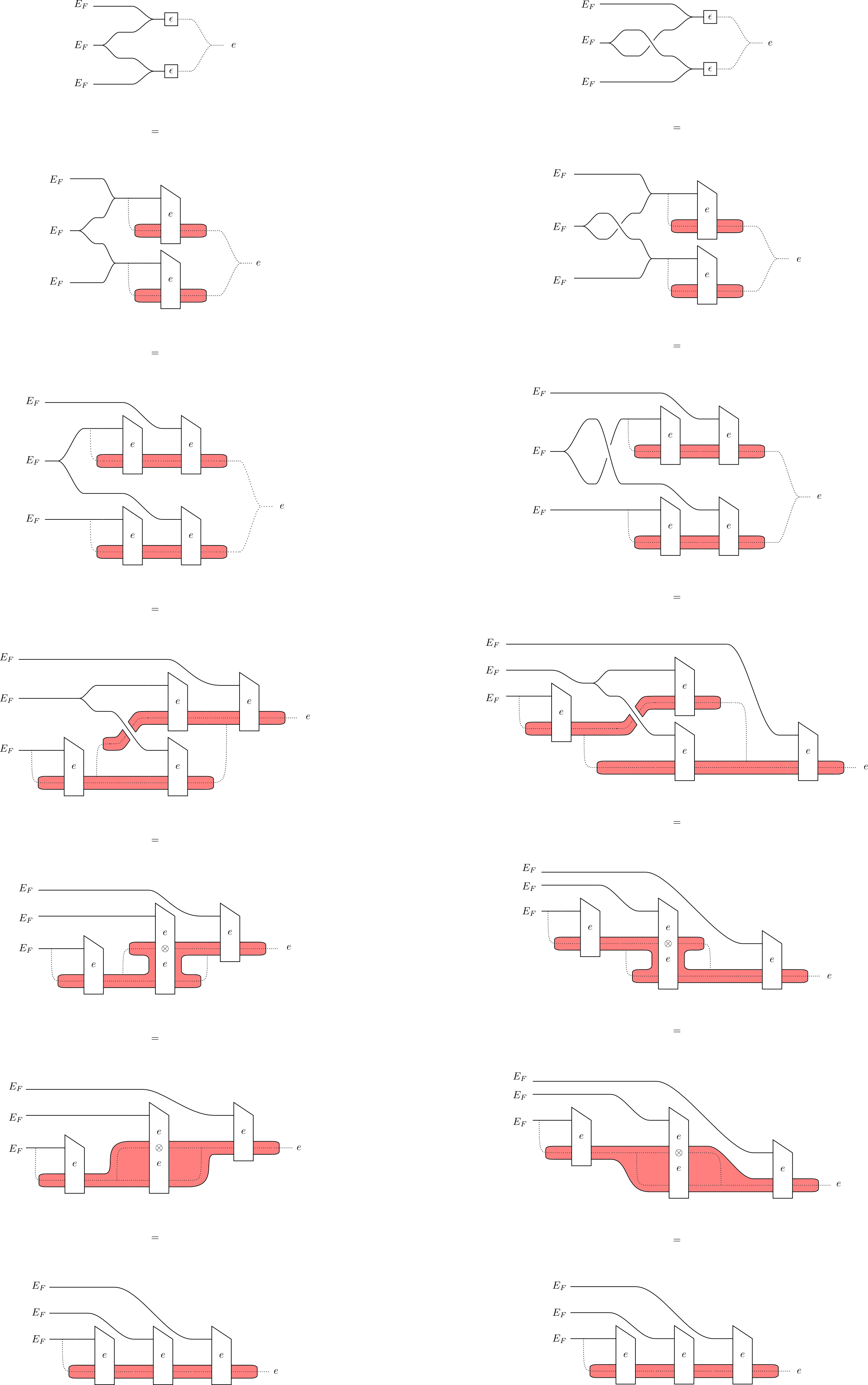}}\end{center}
\caption{Weak counit calculations}
\label{weakcounit23}
\end{figure}

The calculations in figure (\ref{weakunit23}) show that the second and third unit expressions have the following discharged forms:
\begin{center}\resizebox{450pt}{!}{\includegraphics{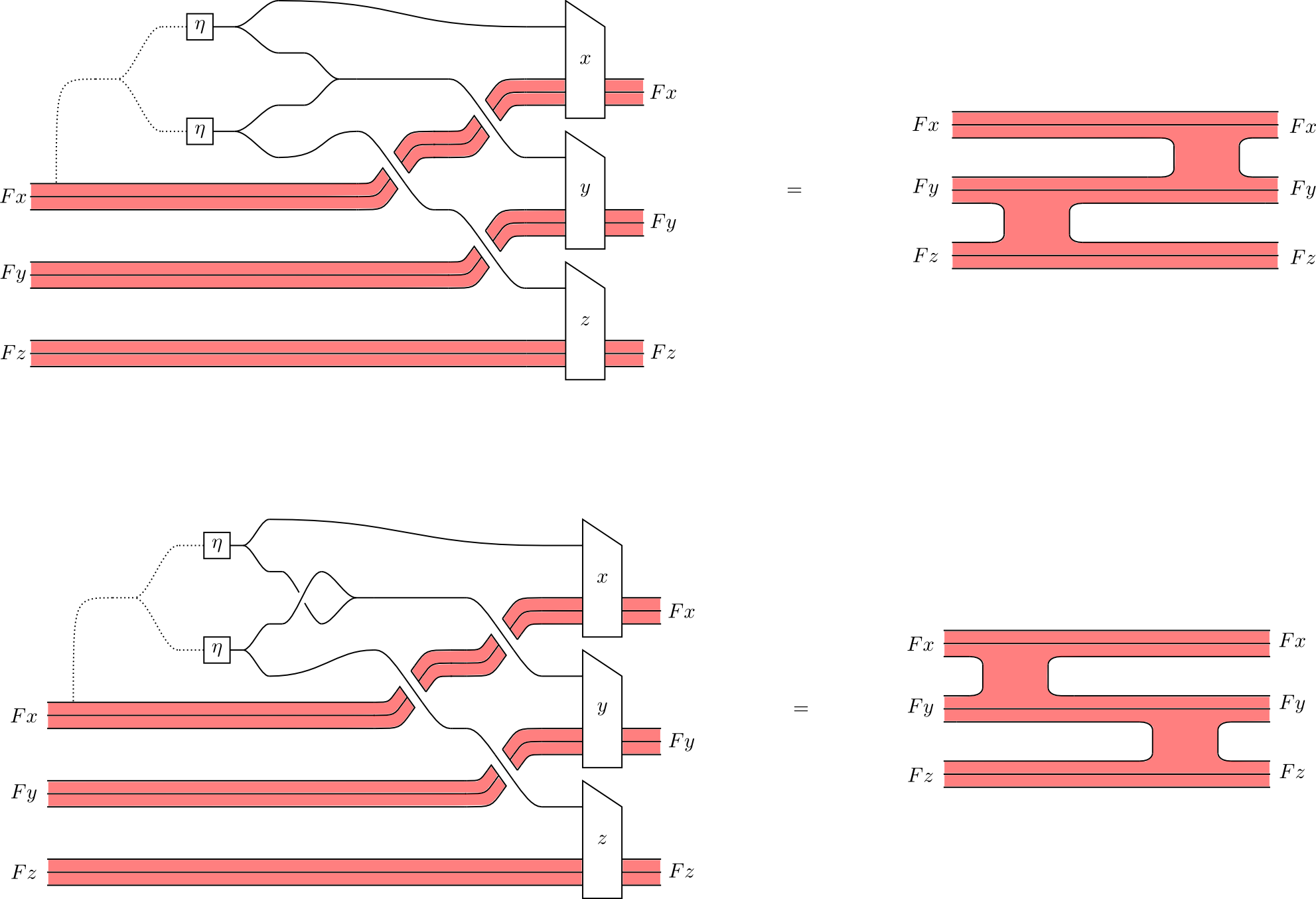}}\end{center}
For these unit axioms, we see that it suffices to assume that $F$ is Frobenius. As for the counit axioms, the discharged form of the first of these is easily calculated:
\begin{center}\resizebox{450pt}{!}{\includegraphics{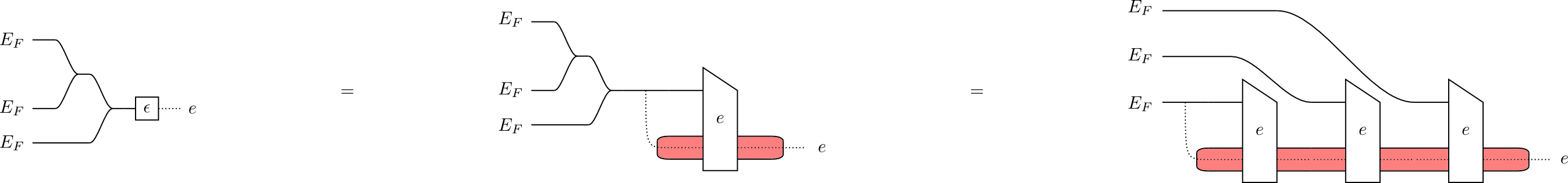}}\end{center} The discharged forms of the second and third counit expression are computed
in figure (\ref{weakcounit23}) and are, of course the same. Examing this figure shows that the counit axioms follow merely from $F$ being both monoidal and comonoidal,
without requiring Frobenius or separable.

This assymmetry (between unit and counit axioms) results from defining $E_F$ using an end, had we instead used a coend, the situation would be reversed.

As for the antipode axioms, we can also consider the pair of strong antipode axioms or the trio of weak antipode axioms.
The strong antipode axioms request the following two equations:
\begin{center}\resizebox{450pt}{!}{\includegraphics{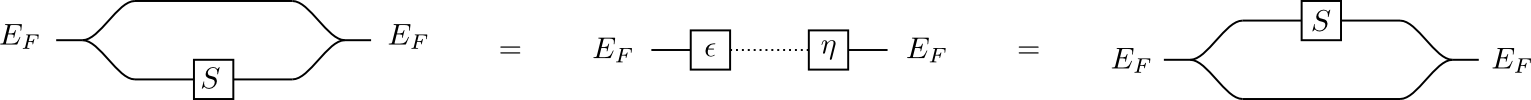}}\end{center} 
On the other hand, the weak antipode axioms request the following three equations:
\begin{center}\resizebox{350pt}{!}{\includegraphics{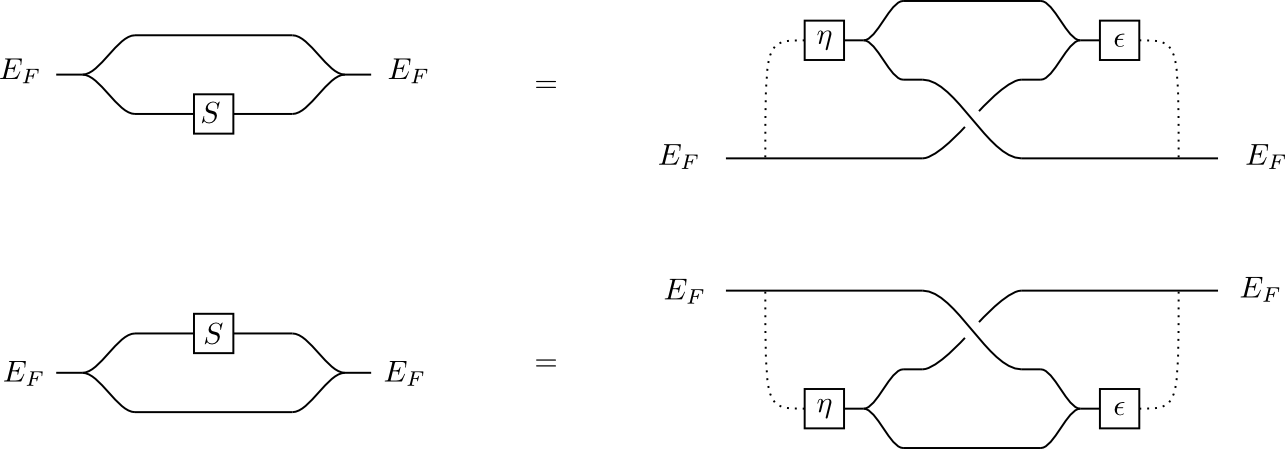}}\end{center} 
\begin{center}\resizebox{350pt}{!}{\includegraphics{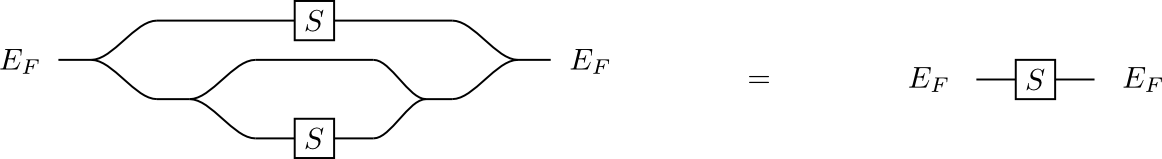}}\end{center} 

\begin{figure}[ht]
\begin{center}\resizebox{!}{550pt}{\includegraphics{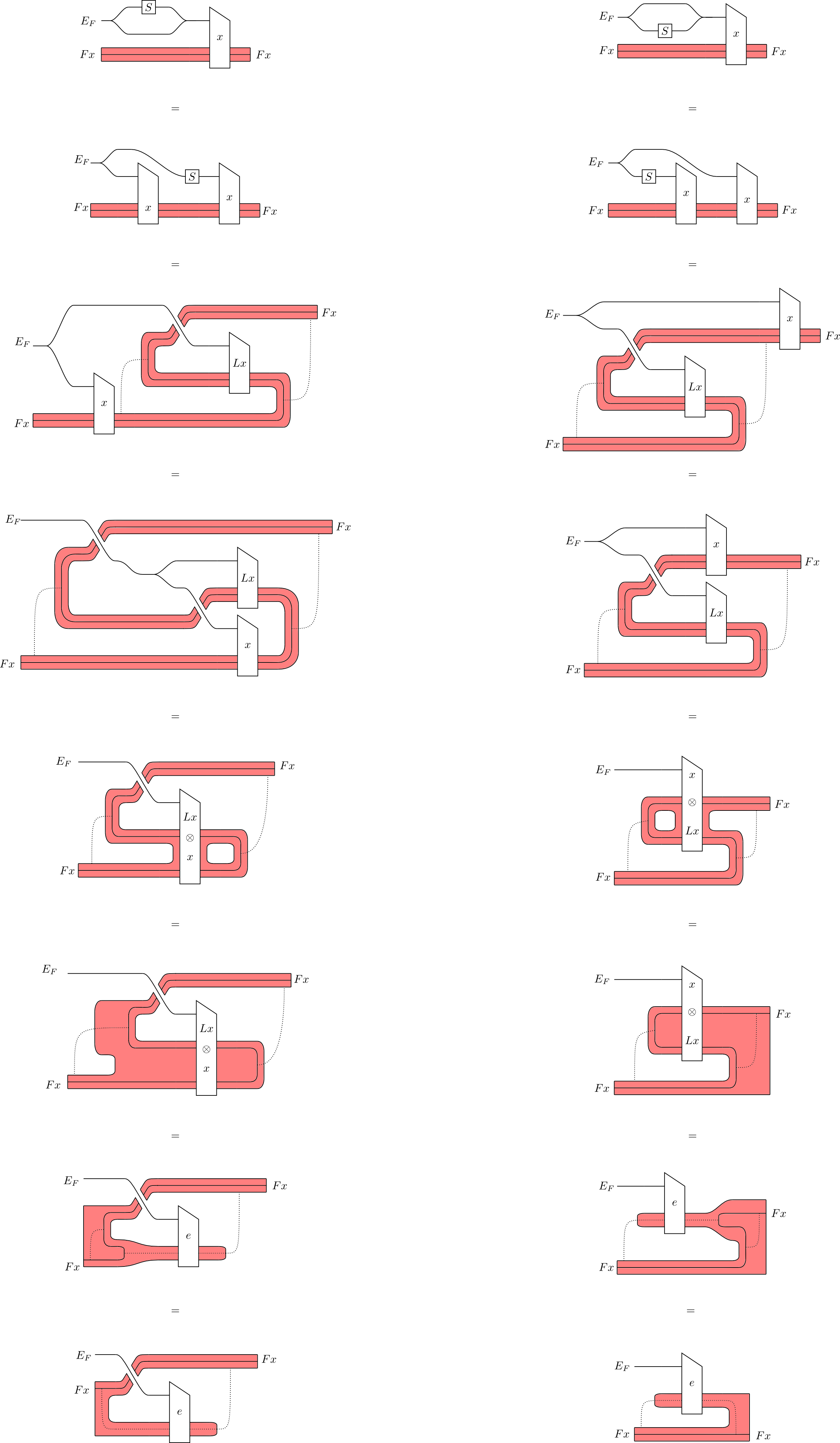}}\end{center}
\caption{Calculations of $S \star E_F$ and $E_F \star S$}
\label{antipodeconvolutions}
\end{figure}

Both sets of axioms involve the two convolutions of the antipode with the (compositional) identity, and so we calculate these two quantities explicitly.
The pair of calculations in figure (\ref{antipodeconvolutions}) show that the discharged forms of $S \star E_F$ and $E_F \star S$ are the following:
\begin{center}\resizebox{450pt}{!}{\includegraphics{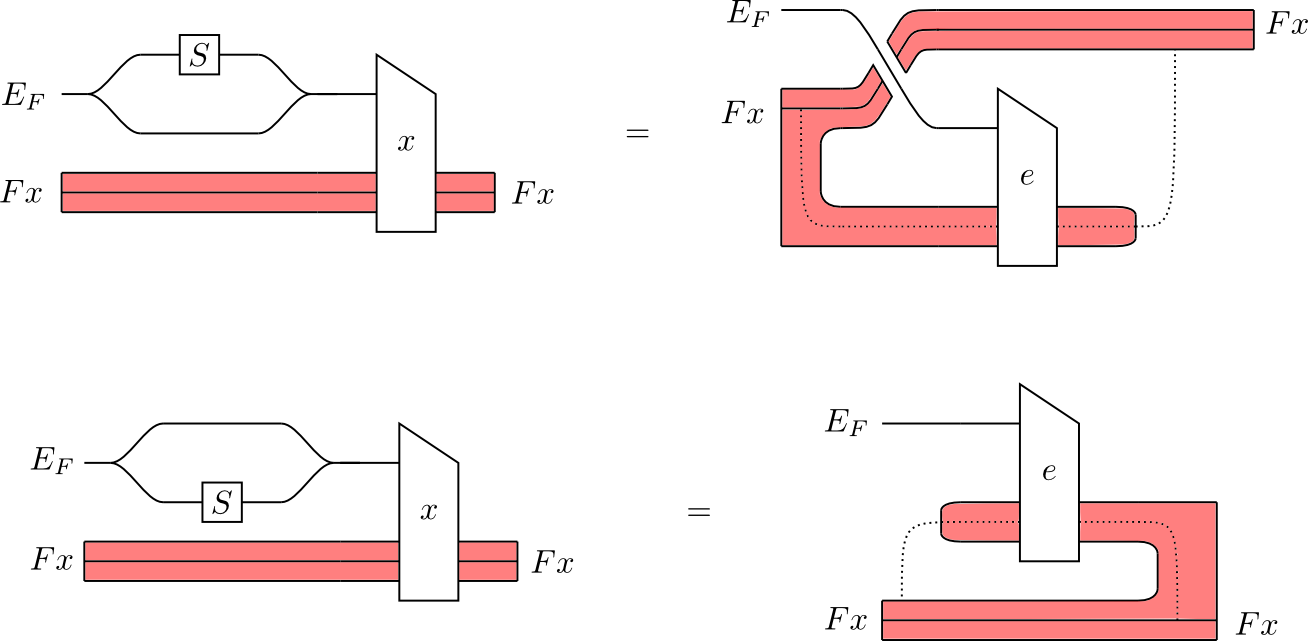}}\end{center}

In the (strong) case, both of these convolutions are supposed to equal the composite $E_F \oto{\epsilon} e \oto{\eta} E_F$, 
the discharged form of which we compute:
\begin{center}\resizebox{350pt}{!}{\includegraphics{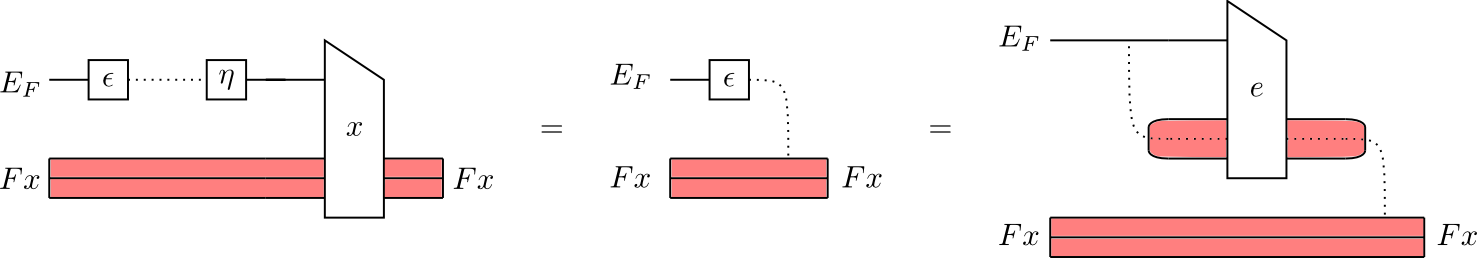}}\end{center}

\begin{figure}[ht]
\begin{center}\resizebox{!}{550pt}{\includegraphics{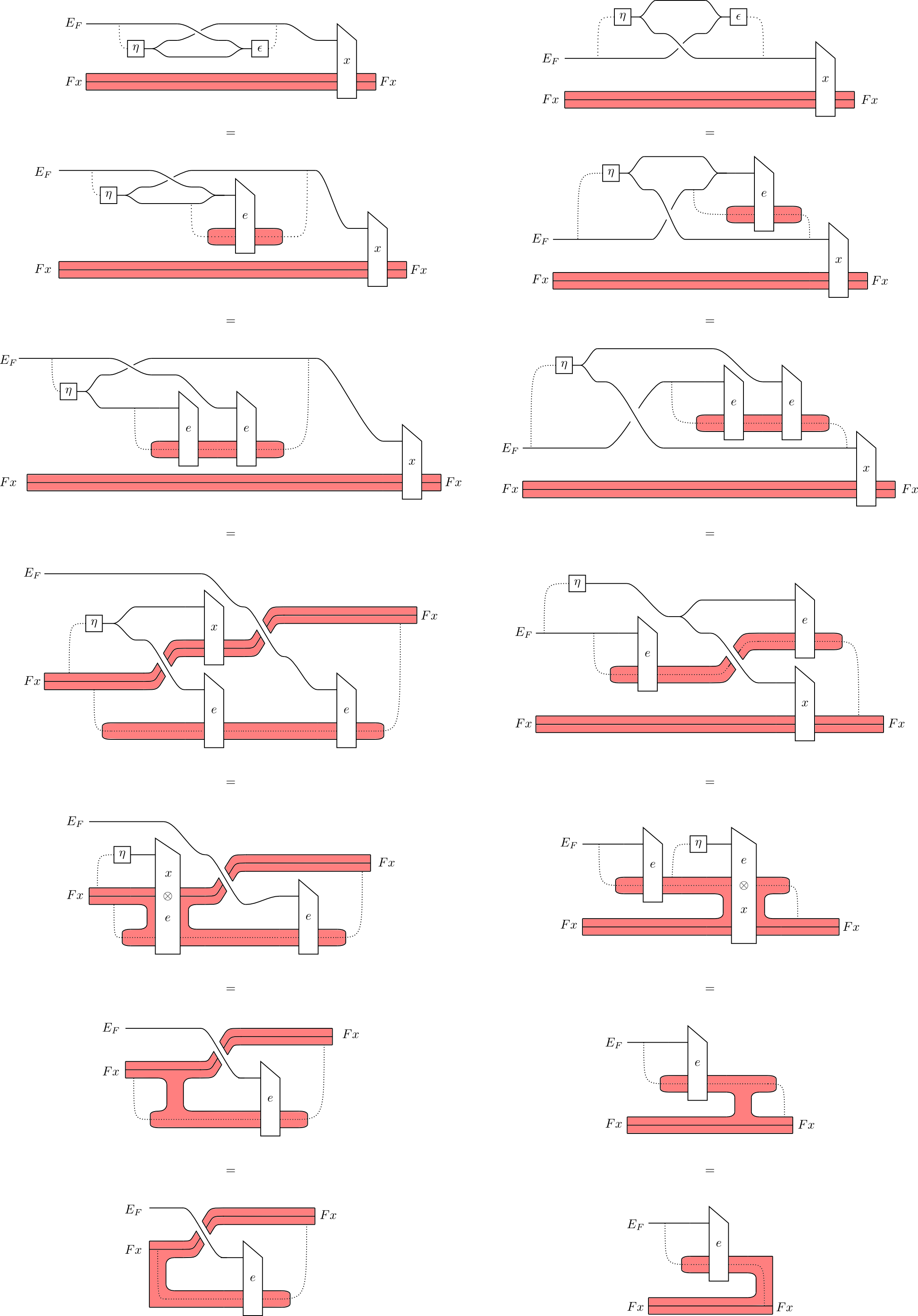}}\end{center}
\caption{``Source'' and ``Target'' maps.}
\label{ants}
\end{figure}

Therefore, we see that, in the usual Hopf algebra case, it suffices to take $Fe \too e \too Fe$ equal to the identity.
For the weak case, these convolutions are instead set equal to the expressions which are computed in figure (\ref{ants}).
So we see that for these two axioms it suffices to take $F$ to be monoidal and comonoidal.

\begin{figure}[ht]
\begin{center}\resizebox{!}{550pt}{\includegraphics{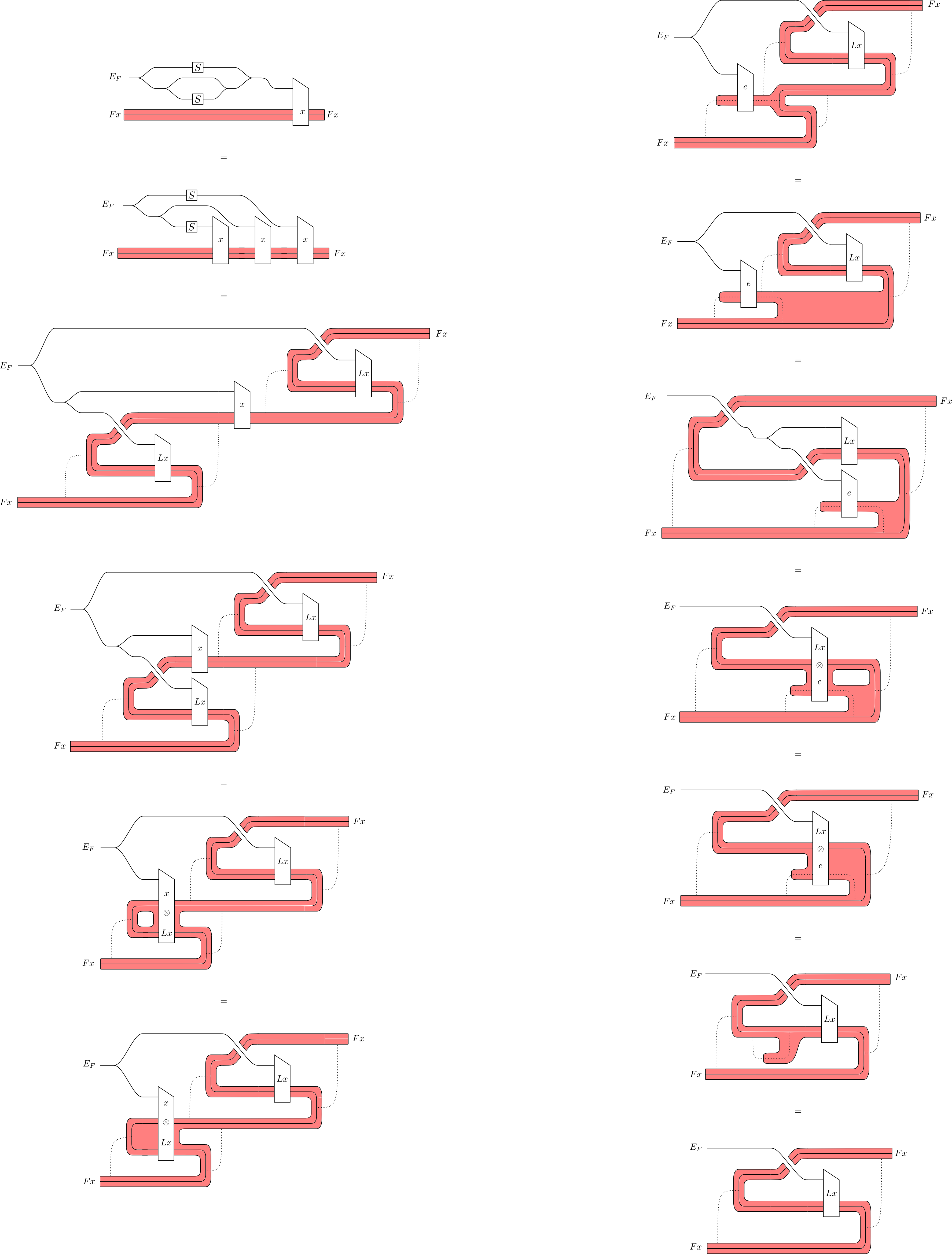}}\end{center}
\caption{The calculation showing $S \star E_F \star S = S$}
\label{triples}
\end{figure}

There is one additional antipode axiom which is imposed for a weak Hopf algebra, namely, that the convolution $S \star E_F \star S$ should equal $S$. 
For this, we compute the left hand side in figure (\ref{triples}), the last diagram of which is the definition of $S$, as desired.



An updated version of this paper will treat reconstruction of braided weak bialgebras and braided weak Hopf algebras when $A$ is known to be
a braided category.

\end{document}